\documentclass[12pt]{article}
\usepackage[cp1251]{inputenc}
\usepackage[T2A]{fontenc}
\usepackage{textcomp}
\usepackage{amsmath, amssymb}

\usepackage[left=3cm,right=3cm,
    top=3cm,bottom=3cm,bindingoffset=0cm]{geometry}
 
\begin{document}

\begin{center}
    \Large {\bf From Rota-Baxter operators to quasitriangular Lie bialgebra structures on $gl_2(\mathbb C)$.}
\end{center}

\begin{center}
{\bf М.\,Goncharov}
\end{center}
\footnotetext{*) This research is supported by Russian Science Foundation (project 21-11-00286),
https://rscf.ru/en/project/21-11-00286/.}

\begin{abstract}
    In the paper we describe structures of quasitriangular Lie bialgebra on $gl_2(\mathbb C)$ using the classification of Rota-Baxter operators of nonzero weight on $gl_2(\mathbb C)$.\bigskip

\end{abstract}

\section{Introduction}

Let $A$ be an arbitrary algebra over a field $F$, $\lambda\in F$. A
map $R:A\mapsto A$ is called a Rota-Baxter operator of weight
$\lambda$ if for all $x,y\in A$:
\begin{equation}\label{e1}
R(x)R(y)=R(R(x)y+xR(y)+\lambda xy).
\end{equation}

Rota-Baxter operators for associative algebras first appear in the
paper by G. Baxter as a tool for studying integral operators that
appear in the theory of probability and mathematical statistics
\cite{Baxter}. For a long period of time, Rota Baxter operators were intensively studied in combinatorics and probability theory mainly. For basic
results and the main properties of Rota---Baxter algebras, see
\cite{Guo}.

Independently, in 80-th Rota-Baxter  operators of weight 0 on Lie algebras
naturally appear in the papers of A.A. Belavin, V.G. Drinfeld \cite{BD}
and M.A. Semenov-Tyan-Shanskii \cite{STS} while studying solutions of the classical Yang-Baxter equation. It was mentioned, that for any quadratic Lie algebra $(L,\omega)$, the standard technique of multilinear algebra gives the one-to-one correspondence between skew-symmetric solutions of the classical Yang-Baxter equation on $L$ and Rota-Baxter operators $R:L\mapsto L$ of weight 0, satisfying $R^*=-R$ ($R^*$ is the adjoint to $R$ operator with respect to the form $\omega)$.

Skew-symmetric solutions of the classical Yang-Baxter equation on a Lie algebra $L$ induce on $L$ the structure of a (triangular) Lie bialgebra. 

If we consider non-skew-symmetric solutions of CYBE, then they induce a structure of a (quasitriangular) Lie bialgebras if and only if the symmetric part of the solution is $L$-invariant. In \cite{GME}, the connection between solutions of this type and Rota-Baxter operators of weight 1 was studied. It was shown that for a finite-dimensional simple Lie algebra $L$ over a field of characteristic 0 there is one-to-one correspondence between non-skew-symmetric solutions of CYBE with $L$-invariant symmetric part and Rota-Baxter operators of weight 1, satisfying 
\begin{equation}\label{LSS}
    L+L^*+\mathrm{id}=0.
\end{equation} 
If $L$ is not simple, then the connection is not straightforward (see Statement 1 and Theorem 1 below). In \cite{Fact}, the class of quasitriangular Lie bialgebras that correspond to Rota-Baxter operators satisfying \eqref{LSS} was introduced and studied.

It is worth noting that  for many varieties of algebras (associative, Jordan, alternative ect.) all structures of corresponding bialgebras on semisimple finite-dimensional  algebras from these varieties are triangular (since they are unital, see, for example, \cite{Zhel} for Jordan algebras). This means that in these varieties, Rota-Baxter operators satisfying \eqref{LSS} do not seem to be interesting (it is known that there are no Rota-Baxter operators of weight 1 on $M_n(F)$ satisfying \eqref{LSS}).

There is a standard method for classification of skew-symmetric solutions of a classical Yang-Baxter equation on a given algebra $A$ (of arbitrary variety): it is known that these solutions are in one-to-one correspondence with pairs $(B,\omega)$, where $B$ is a subalgebra in $A$ and $\omega$ is a symplectic form on $B$ (see \cite{BD}). But there is no standard method for classification of structures of quasitriangular bialgebras that may be used for an arbitrary variety of algebras  (where solutions of this type are interesting). Here we can mention the paper \cite{Stolin}, where the classification was made under the assumption that $L$ is a simple complex Lie algebra and involved some information about root systems on $L$. If $M$ is a Malcev algebra, the classification of quasitriangular structures on $M$ was made in \cite{GME2} by considering some specific information concerning the classical limit $M\oplus M^*$.

In recent years, there have been many papers where Rota-Baxter operators on important classes of algebras were described (\cite{sl2-0},\cite{sl2},\cite{GonGub}, ect.). Usually, the description is made up to the action of the group of automorphisms. The natural question here is,  if we can use these results to describe the structures of a quasitriangular bialgebra on these classes of Lie (Malcev, ect.) algebras?  Unfortunately, if the description of Rota-Baxter operators was made "up to an automorphism", then we can't use it directly since conjugate operators not necessary give conjugate tensors. In the current paper we consider the question for general linear Lie  algebra of order 2 $gl_2(\mathbb C)$ over the field of complex numbers. We use the description of the Rota-Baxter operators  on $gl_2(\mathbb C)$ obtained in \cite{GME1} to classify (up to the action of $Aut(gl_2(\mathbb C))$) structures of quasitriangular Lie bialgebras on $gl_2(\mathbb C)$.

\section{Motivation and preliminary results.}

All vector spaces are assumed to be over a field of complex numbers $\mathbb C$. Given a vector space $V$  over $\mathbb C$, denote by
$V\otimes V$ it's tensor product over $\mathbb C$. Define the linear mapping
$\tau$ on $V\otimes V$ by  $\tau(\sum\limits_ia_i\otimes
b_i)=\sum\limits_ib_i\otimes a_i$.

Let $L$ be a Lie algebra with a product $[\cdot,\cdot]$. 

A Lie algebra $L$ acts on
$L^{\otimes n}$ by
$$
[x_1\otimes x_2\otimes \ldots \otimes x_n,y]=\sum\limits_i
x_1\otimes\ldots\otimes [x_i,y]\otimes \ldots\otimes x_n$$
for all $x_i,y\in L$.

{\bf Definition.} An element $r\in L^{\otimes n}$ is called
$L$-invariant (or $ad$-invariant) if $[r,y]=0$ for all $y\in L$.

{\bf Definition.} We say that an element $r=\sum_i a_i\otimes b_i\in L\otimes L$ is a solution of the classical Yang-Baxter equation  (CYBE) on $L$ if

\begin{equation}\label{lieYB}
\sum\limits_{ij}[a_i,a_j]\otimes b_i\otimes
b_j-a_i\otimes[a_j,b_i]\otimes b_j+a_i\otimes a_j\otimes
[b_i,b_j]=0.
\end{equation}

A solution $r\in L\otimes L$ of CYBE is called skew-symmetric, if $\tau(r)=-r$. 

{\bf Remark.} Let $r=\sum a_i\otimes b_i \in L\otimes L$ be a solution of CYBE and $\varphi\in Aut(L)$, then $r_1=\sum \varphi(a_i)\otimes \varphi(b_i)$ is also a solution of CYBE on $L$. In this case we will say that tensors $r$ and $r_1$ are conjugate. Moreover, if $r$ is skew-symmetric (the symmetric part of $r$ is $L$-invariant), then so is $r_1$ (so is the symmetric part of $r_1$). Therefore, it is possible to make a description of CYBE up to the action of $Aut(L)$.

{\bf Definition.} A bilinear form
$\omega$ on a Lie algebra $L$ is called invariant if
$\omega([a,b],c)=\omega(a,[b,c])$ for all $a,b,c\in L$.

{\bf Definition.} Let $L$ be a Lie algebra and $\omega$ be an
invariant non-degenerate form on $L$. Then the pair $(L,\omega)$ is
called a quadratic Lie algebra.

Given  a quadratic Lie algebra $(L,\omega)$, for every element $r=\sum\limits_i a_i\otimes b_i\in L\otimes L$ we may define a linear map $R: L\mapsto L$ as
\begin{equation}\label{rb2}
    R(a)=\sum\limits_i \omega (a_i,a)b_i,
\end{equation}
 $a\in L$. By $R^*$ define the adjoint map with respect to the form $\omega$: $$\omega(R(a),b)=\omega(a,R^*(b))$$
 for all $a,b\in L$.
 
 It is known that an element $r\in L\otimes L$ is a skew-symmetric solution of CYBE on a quadratic Lie algebra $L$ if and only if the corresponding map $R:L\mapsto L$ is a Rota-Baxter operator of weight 0 satisfying $R+R^*=0$ \cite{STS}.

{\bf Definition \cite{Drinf}.}  Let $L$ be a Lie algebra with a comultiplication
$\delta: L\mapsto L\otimes L$. The pair $(L,\delta)$ is called a Lie bialgebra if and
only if $(L,\delta)$ is a Lie coalgebra and $\delta$ is a 1-cocycle,
i.e., it satisfies
\begin{multline*}
\delta([a,b])= [\delta(a),b]+[a,\delta(b)]=\\
  =\sum([a_{(1)},b]\otimes
a_{(2)}+a_{(1)}\otimes[a_{(2)},b])+\sum([a,b_{(1)}]\otimes
b_{(2)}+b_{(1)}\otimes[a,b_{(2)}]).    
\end{multline*}

for all $a,b\in L$.

There is an important type of Lie bialgebras: let $L$ be a Lie algebra and
$r=\sum\limits_ia_i\otimes b_i\in L\otimes L$. Define a
comultiplication $\delta_r$ on $L$ by
$$
\delta_r(a)=[r,a]=\sum\limits_i[a_i,a]\otimes b_i+a_i\otimes[b_i,a]
$$
for all $a\in L$. It is easy to see that $\delta_r$ is a 1-cocycle.

 The dual algebra $L^*$ of the coalgebra
$(L,\delta_r)$ is skew-symmetric if and only if $r+\tau(r)$ is
$L$-invariant. Also, $L^*$  satisfies the Jacobi identity if and
only if an element $C_L(r)$, defined as
$$
C_L(r)=\sum\limits_{ij}[a_i,a_j]\otimes b_i\otimes
b_j-a_i\otimes[a_j,b_i]\otimes b_j+a_i\otimes a_j\otimes
[b_i,b_j],
$$
is $L$-invariant. 

If $\tau(r)=-r$ and $r$ is a solution of CYBE, then $(L,\delta_r)$ is a Lie bialgebra called {\it triangular} Lie bialgebra. If $r+\tau(r)$ is a nonzero $L$-invariant element and $r$ is a solution of CYBE, then $(L,\delta_r)$ is called a {\it quasitriangular} Lie bialgebra. Triangular and quasitriangular Lie bialgebras play an important role since they lead to solutions of the quantum Yang-Baxter equation.

{\bf Remark.} Let $r\in L\otimes L$, $\varphi\in Aut(L)$ and $r_1=(\varphi\otimes \varphi)r$. It is easy to see that bialgebras $(L,\delta_r)$ and $(L,\delta_{r_1})$ are isomorphic. Therefore, the description of CYBE up to action of $Aut(L)$ agrees with the classification of structures of quasitriangular Lie bialgebra structures on $L$.

Let  $(L,\omega)$ be a quadratic Lie algebra over a field $F$ and $r=\sum\limits_i a_i\otimes b_i\in L\otimes L$. Let $R$ be a linear map defined as in \eqref{rb2} and
$R^*$ be the adjoint map with respect to the form $\omega$. In general, if $L$ is not a simple complex Lie algebra, $r$ is a solution of CYBE and $r+\tau(r)$ is $L$-invariant, then $R$ is not necessary a Rota-Baxter operator (see \cite{GME} for details and examples).  In what follows, we will need the following Statement 1 and Theorem 1 from \cite{GME}. These results give the connection between solutions of CYBE and Rota-Baxter operators on $L$.
 
 {\bf Statement 1.} The symmetric part $r+ \tau(r)$ of $r$ is $L$-invariant if and only if for all $a,b\in L$ 
 \begin{equation}\label{rb3}
    R([a,b])+R^*([a,b])=[R(a)+R^*(a),b].
\end{equation}

{\bf Theorem 1.} 1. If $r$ is a solution of the classical Yang-Baxter equation on $L$ then $R$ is a Rota-Baxter operator of weight 1 if and only if for all $a,b\in L$:
\begin{equation}\label{rb1}
[R(a),b]+[R^*(a),b] + [a,b]\in ker(R).    
\end{equation} 
2. Let $R:L\mapsto L$ be a Rota-Baxter operator of weight 1 and let $r\in L\otimes L$ be the corresponding tensor defined as in \eqref{rb2}. Then $r$ is a solution of the classical Yang-Baxter equation if and only if $R$ satisfies  \eqref{rb1}.

\section{Main result.}

Let $M_2(\mathbb C)$ be the matrix algebra of order 2 over $\mathbb C$ with the multiplication $xy$. The multiplication in the general linear algebra $gl_2(\mathbb C)$ we will denote by $[\cdot,\cdot]$: 
$$
[x,y]=xy-yx,
$$
$x,y\in gl_2(\mathbb C)$.
We will describe all solutions of the classical Yang-Baxter equation with $ad$-invariant even part on the general linear algebra of order 2 $gl_2(\mathbb C)$. Note, that $gl_2(\mathbb C)$ contains a nontrivial center $Z$ spanned by the identity matrix $\mathrm E$ and is not a semisimple Lie algebra. We will also consider $sl_2(\mathbb C)$ as a Lie subalgebra in $gl_2(\mathbb C)$. Then $gl_2(\mathbb C)=\mathbb C \mathrm E\oplus sl_2(\mathbb C)$ and for any $\varphi\in Aut(gl_2(\mathbb C))$: $\varphi(\mathrm E)=\theta \mathrm E$, $\varphi(sl_2(\mathbb C))=sl_2(\mathbb C)$.

We will consider $gl_2(\mathbb C)$ as a quadratic Lie algebra with the trace form $\omega$ defined as
$$
\omega(x,y)=\mathrm{tr}(xy)
$$

{\bf Theorem 2.} Let $r\in gl_n(\mathbb C)\otimes gl_n(\mathbb C)$ be a solution of CYBE with $gl_n(\mathbb C)$-invariant nonzero even part $r+\tau(r)\neq 0$. Then the corresponding map $R$ defined by \eqref{rb2}  is a Rota-Baxter operator of a nonzero weight $\lambda$. Moreover, the adjoint map $R^*$ is also a Rota-Baxter operator of the same weight and for any $x\in sl_n(\mathbb C)$ we have
$$
R(x)+R^*(x)+\lambda id=0.
$$

{\bf Proof.}
For any $\lambda\in \mathbb C$ consider a map $\theta_{\lambda}:gl_n(\mathbb C)\mapsto gl_n(\mathbb C)$: for any $x\in gl_n(\mathbb C)$ put
$$
\theta_{\lambda}(x)=R(x)+R^*(x)+\lambda x.
$$

Consider a set 
$$
I_{\lambda}=\{\theta_{\lambda}(x)|\ x\in [gl_n(\mathbb C),gl_n(\mathbb C)]\}.
$$ 

In \cite{GME} it was proved that $I_{\lambda}$ is an ideal in $gl_n(\mathbb C)$ for any $\lambda$. Moreover, $I_{\lambda}\subset [gl_n(\mathbb C),gl_n(\mathbb C)]=sl_n(\mathbb C)$ (and as consequence, $sl_n(\mathbb C)$ is $\theta_{\lambda}$-invariant). Since $sl_n(\mathbb C)$ is simple, we have two possibilities: $I_{\lambda}=0$ or $I_{\lambda}=sl_n(\mathbb C)$. We want to prove that there exists a unique $\alpha\in \mathbb C$ such that $I_{\alpha}=0$. The uniqueness is quite straightforward: if $I_{\alpha_1}=I_{\alpha_2}=0$, then for any $x\in sl_n(\mathbb C)$: 
$$
R(x)+R^*(x)+\alpha_1 x=R(x)+R^*(x)+\alpha_2 x
$$
that is not possible if $\alpha_1\neq \alpha_2$. 

Take arbitrary $\lambda\in \mathbb C$. Since $r+\tau(r)$ is $gl_n(\mathbb C)$-invariant, then $\theta_{\lambda}$ satisfy 
$$
\theta_{\lambda}([x,y])=[\theta_{\lambda}(x),y]
$$
for all $x,y\in gl_n(\mathbb C)$. In other words, $\theta_{\lambda}$ belongs to the centraliser of $gl_n(\mathbb C)$. Since $sl_n(\mathbb C)$ is a simple complex Lie algebra, the restriction of $\theta_{\lambda}$ to $sl_n(\mathbb C)$ is equal to $\gamma id$ for some $\gamma\in \mathbb C$. It means that $I_{\lambda-\gamma}=\theta_{\lambda-\gamma}(sl_n(\mathbb C))=0$.

Take the scalar $\lambda\in \mathbb C$ such that $I_{\lambda}=0$. Obviously, $I_{\lambda}$ is $R$-invariant. Now we can use Theorem 2 form \cite{GME} to get that $R$ and $R^*$ are Rota-Baxter operators of weight $\lambda$ on the quotient algebra $gl_n(\mathbb C)/I_{\lambda}=gl_n(\mathbb C)$. If $\lambda=0$, then it is known that the tensor $r$ is skew-symmetric, so $\lambda\neq 0$. The theorem is proved.

{\bf Remark.} Apart from the case of simple complex Lie algebras, here we can't say that 
$$
R+R^*+\lambda id=0.
$$
 In Theorem 2 we proved that for all $x\in sl_n(\mathbb C)$: $\theta_{\lambda}(x)=R(x)+R^*(x)+\lambda x=0$. But in general $\theta_{\lambda}(\mathrm E)\neq 0$.

Note that the map $R$ from Theorem 2 should satisfy \eqref{rb3} and \eqref{rb1}. By statement 1 and Theorem 1 the inverse is also true: if $R$ is a Rota-Baxter operator of weight 1 satisfying \eqref{rb3} and \eqref{rb1}, then the corresponding tensor $r$ is a non-skew-symmetric solution of the CYBE with $gl_n(\mathbb C)$-invariant even part $r+\tau(r)$. Thus, we need to classify such operators. 

In \cite{GME1} the description of Rota-Baxter operators of weight 1 on $gl_2(\mathbb C)$ was given. The description was given up to the conjugation with automorphisms from $Aut(gl_2(\mathbb C))$.

{\bf Theorem 3.} Let $R$ be a Rota-Baxter operator of weight 1 on $gl_2(\mathbb C)$. Then, up to conjugation with automorphisms of $Aut(gl_2(\mathbb C))$, $R$ is equal to one of the following

\begin{gather*}
1.\  R(\mathrm E)=\lambda \mathrm E+e_{12},\ R(h)=R(e_{12})=R(e_{21})=0;\\
2.\ R(\mathrm E)=\lambda \mathrm E+e_{12},\ R(h)=-h,\  R(e_{12})=-e_{12},\  R(e_{21})=-e_{21};\\
3.\  R(\mathrm{E})=\lambda \mathrm{E}+h,\ R(h)=0,\ R(e_{12})=R(e_{21})=0,\ \lambda\in\mathbb C ;\\
4.\  R(\mathrm{E})=\lambda \mathrm{E}+h,\ R(h)=-h,\ R(e_{12})=-e_{12},\ R(e_{21})=-e_{21},\ \lambda\in\mathbb C;\\
 5.\  R(\mathrm{E})=\lambda \mathrm{E}+h,\ R(h)=\alpha_1\mathrm{E}+\alpha_2h,\ R(e_{12})=-e_{12},\ R(e_{21})=0,\ \lambda\in\mathbb C;\\
6.\  R(\mathrm{E})=\lambda \mathrm{E},\ R(h)=0,\ R(e_{12})=-e_{12}+th;\ R(e_{21})=0,\ t\in\{0,1\};\\
7.\  R(\mathrm{E})=\lambda \mathrm{E},\ R(h)=0,\ R(e_{12})=-e_{12}+th+E;\ R(e_{21})=0,\ t\in\{0,1\};\\
 8.\  R(\mathrm{E})=\lambda \mathrm{E},\ R(h)=E,\ R(e_{12})=-e_{12}+h+\alpha \mathrm{E};\ R(e_{21})=0,\ \alpha\in\mathbb C;\\
9.\  R(\mathrm{E})=\lambda \mathrm{E},\ R(h)=E,\ R(e_{12})=-e_{12}+\mathrm{E};\ R(e_{21})=0;\\
10.\  R(\mathrm{E})=\lambda \mathrm{E},\ R(h)=th,\ R(e_{21})=0,\ R(e_{12})=-e_{12},\ t\in \mathbb C,\ t\neq 0;\\
11.\  R(\mathrm{E})=\lambda \mathrm{E},\ R(h)=th+\mathrm{E},\ R(e_{21})=0,\ R(e_{12})=-e_{12},\ t\in \mathbb C,\ t\neq 0;\\
 12.\  R(\mathrm{E})=\lambda \mathrm{E},\ R(h)=-h+\alpha \mathrm{E},\ R(e_{21})=\mathrm{E},\ R(e_{12})=-e_{12},\ \alpha\in \mathbb C;\\
   13.\ R(\mathrm{E})=\lambda\mathrm{E}, R(h)=th,\ R(e_{12})=te_{12},\ R(e_{21})=te_{21},\ t\in\{0,-1\}.
\end{gather*}

{\bf Remark 2.} Unfortunately, we can't use this result to describe all solutions of the modified classical Yang-Baxter equation on $gl_2(\mathbb C)$. Indeed, let $R$ be a Rota-Baxter operator on $gl_2(\mathbb C)$, $r=\sum a_i\otimes b_i$ be the corresponding tensor and  $\varphi$  be an automorphism from $Aut(gl_2(\mathbb C))$. Consider $R_1=\varphi^{-1}\circ R\circ\varphi$. Then the corresponding to $R_1$  tensor is the following:
$$r_1=\sum \varphi^* (a_i)\otimes \varphi^{-1}(b_i).$$
In other words, tensors $r$ and $r_1$ are not necessarily conjugate by an automorphism from $Aut(gl_2(\mathbb C))$. There may be a situation when $r$ is a solution of CYBE while $r_1$ is not a solution.

Nevertheless, we have the following 

{\bf Proposition 1.} If $\varphi$ is an automorphism of $M_2(\mathbb C)$ (as an associative algebra), then $\varphi$ is an "orthogonal" automorphism of $gl_2(\mathbb C)$, that is $\varphi^*=\varphi^{-1}$. This means that if $R:gl_2(\mathbb C)\mapsto gl_2(\mathbb C)$ is a linear map and $R_1=\varphi^{-1}\circ R\circ \varphi$, then corresponding tensors $r$ and $r_1$ (to $R$ and $R_1$ respectively) are conjugate:
$$
r_1=(\varphi^{-1}\otimes \varphi^{-1})r.
$$
{\bf Proof.}
Indeed, if $\varphi\in Aut(M_2(\mathbb C))$, then 
$$\omega(\varphi(x),\varphi(y))=\mathrm{tr}(\varphi(x)\varphi(y))=\mathrm{tr}(\varphi(xy))=\mathrm{tr}(xy)=\omega(x,y).$$
That means that $\varphi^*=\varphi^{-1}$ and the proposition is proved.

For any $\theta\in \mathbb C$, $\theta\neq 0$, we can define an automorphism $\psi_{\theta}$ of $gl_2(\mathbb C)$ as follows:
\begin{equation}\label{aut}
\psi_{\theta}(\mathrm{E})=\theta \mathrm{E},\ \psi_{\theta}(a)=a 
\end{equation}
for any $a$ satisfying $\mathrm{tr}(a)=0$. It is clear that $\psi_{\theta}\circ\varphi=\varphi\circ\psi_{\theta}$ for all $\varphi\in Aut(gl_2(\mathbb C))$.

{\bf Proposition 2.}
For any $\varphi\in Aut(gl_2(\mathbb C))$ there are $0\neq\theta\in \mathbb C$ and $\phi\in Aut(M_2(\mathbb C))$ such that 
$$
\varphi=\psi_{\theta}\circ \phi=\phi\circ \psi_{\theta}.
$$
{\bf Proof.}
Let $\varphi\in Aut(gl_2(\mathbb C))$ and $\varphi(\mathrm{E})=\theta\mathrm{E}$. Consider
$$
\phi=\psi_{\theta^{-1}}\circ \varphi.
$$
Then $\phi\in Aut(gl_2(\mathbb C))$ and $\phi(\mathrm(E))=E$. The restriction $\phi|_{sl_2(\mathbb C)}$ is an automorphism of $sl_2(\mathbb C)$. Since all automorphisms of $sl_2(\mathbb C)$ are inner, there is an invertible matrix $u\in M_2(\mathbb C)$ such that $\phi|_{sl_2(\mathbb C)}(a)=uau^{-1}$ for all $a\in sl_2(\mathbb C)$. Therefore, we can conclude that $\phi(a)=uau^{-1}$ for all $a\in gl_2(\mathbb C)$ and $\phi\in Aut(M_2(\mathbb C))$. The proposition is proved.

 Consider the description of Rota-Baxter operators of weight 1 on $gl_2(\mathbb C)$ modulo the action of the group $Aut(M_2(\mathbb C))$. For this, we need to take the representatives $R$ of the orbits from Theorem 3, take $0\neq\theta\in \mathbb C$ and consider the action $\psi_{\theta}^{-1}\circ R\circ \psi_{\theta}$. We get the following

{\bf Theorem 4.} Let $R$ be a Rota-Baxter operator of weight 1 on $gl_2(\mathbb C)$. Then there is a unique  number $i\in \{1,\ldots,13\}$  such that $R=\psi^{-1}\circ R_1\circ \psi$, where $\psi \in Aut(M_2(\mathbb C))$ and $R_1$ belong to the line number $i$ below for some $0\neq\theta\in \mathbb C$ ($\theta$ is not unique)

\begin{gather*}
1.\  R(\mathrm E)=\lambda \mathrm E+\theta e_{12},\ R(h)=R(e_{12})=R(e_{21})=0,\ \theta\neq 0;\\
2.\ R(\mathrm E)=\lambda \mathrm E+\theta e_{12},\ R(h)=-h,\  R(e_{12})=-e_{12},\  R(e_{21})=-e_{21}\ \theta\neq 0;\\
3.\  R(\mathrm{E})=\lambda \mathrm{E}+\theta h,\ R(h)=0,\ R(e_{12})=R(e_{21})=0,\ \lambda\in\mathbb C \ \theta\neq 0;\\
4.\  R(\mathrm{E})=\lambda \mathrm{E}+\theta h,\ R(h)=-h,\ R(e_{12})=-e_{12},\ R(e_{21})=-e_{21},\ \lambda\in\mathbb C\ \theta\neq 0;\\
 5.\  R(\mathrm{E})=\lambda \mathrm{E}+\theta h,\ R(h)=\alpha_1\mathrm{E}+\alpha_2h,\ R(e_{12})=-e_{12},\ R(e_{21})=0,\ \lambda,\alpha_i\in\mathbb C \ \theta\neq 0;\\
6.\  R(\mathrm{E})=\lambda \mathrm{E},\ R(h)=0,\ R(e_{12})=-e_{12}+th;\ R(e_{21})=0,\ t\in\{0,1\};\\
7.\  R(\mathrm{E})=\lambda \mathrm{E},\ R(h)=0,\ R(e_{12})=-e_{12}+th+\theta E;\ R(e_{21})=0,\ t\in\{0,1\}\ \theta\neq 0;\\
 8.\  R(\mathrm{E})=\lambda \mathrm{E},\ R(h)=\theta E,\ R(e_{12})=-e_{12}+h+\alpha \mathrm{E};\ R(e_{21})=0,\ \alpha\in\mathbb C\ \theta\neq 0;\\
9.\  R(\mathrm{E})=\lambda \mathrm{E},\ R(h)=\theta E,\ R(e_{12})=-e_{12}+\theta\mathrm{E};\ R(e_{21})=0\ \theta\neq 0;\\
10.\  R(\mathrm{E})=\lambda \mathrm{E},\ R(h)=th,\ R(e_{21})=0,\ R(e_{12})=-e_{12},\ t\in \mathbb C,\ t\neq 0;\\
11.\  R(\mathrm{E})=\lambda \mathrm{E},\ R(h)=th+\theta \mathrm{E},\ R(e_{21})=0,\ R(e_{12})=-e_{12},\ t\in \mathbb C,\ t\neq 0\ \theta\neq 0;\\
 12.\  R(\mathrm{E})=\lambda \mathrm{E},\ R(h)=-h+\alpha \mathrm{E},\ R(e_{21})=\theta\mathrm{E},\ R(e_{12})=-e_{12},\ \alpha\in \mathbb C\ \theta\neq 0;\\
   13.\ R(\mathrm{E})=\lambda\mathrm{E}, R(h)=th,\ R(e_{12})=te_{12},\ R(e_{21})=te_{21},\ t\in\{0,-1\}.
\end{gather*}

{\bf Remark.} Here, different scalars $\theta$ not necessarily gives us different orbits. For example, it can be shown that in the line 1 it is possible to take $\theta=1$ modulo $Aut(M_2(\mathbb C))$  since the conjugation of $R$ by $\psi_{\theta}$ is equal to the conjugation of $R$ by $\varphi_A$, where $\varphi_A(x)=AxA^{-1}$ for every $x\in M_2(\mathbb C)$ and $A=\begin{pmatrix}\theta & 0\\ 0 & 1\end{pmatrix}$. But for our purposes, it is enough to make such a rough description. 

Maps that lie in the same orbit (modulo $Aut(M_2(\mathbb C)$) from Theorem 4 correspond to isomorphic tensors. Note that a map $R$ satisfies \eqref{rb3} or \eqref{rb1} if and only if $\varphi^{-1}\circ R\circ \varphi$ satisfies the same conditions.  Thus, it is enough to consider one map from every orbit in Theorem 4.

Let $R:gl_2(\mathbb C)\mapsto gl_2(\mathbb C)$ be a Rota-Baxter operator of weight 1 and $r=\sum a_i\otimes b_i\in gl_2(\mathbb C)\otimes gl_2(\mathbb C)$. From Statement 1 it is follows that if $r+\tau(r)$ is $gl_2(\mathbb C)$ invariant, then
\begin{equation}\label{usE}
R(\mathrm{{E}})+R^*(\mathrm{E})=\gamma \mathrm{E}
\end{equation}
for some $\gamma\in \mathbb C$.

{\bf Proposition 3.} Let $R$ be a Rota-Baxter operator of weight 1 satisfying \eqref{usE}. Then $R$ is conjugate by means of $Aut(M_2(\mathbb C))$ to an operator from lines 5 (with $\alpha_1=-\theta$), 6, 10 or 13 in Theorem 4 only.

{\bf Proof.}
Consider $\mathrm E$, $h=e_{11}-e_{22}$, $e_{12}$ and $e_{21}$ as a basis of $gl_2(\mathbb C)$. To check the condition \eqref{usE}, we need to compute $R^*(\mathrm E)$. For this, we need to find $v\in gl_2(\mathbb C)$ such that $\mathrm{tr}(\mathrm{E}R(v))\neq 0$. 

Suppose that $R$ lies in a orbit from line 1. Then $\mathrm{tr}(\mathrm ER(v))\neq 0$ if and only if $v=\gamma \mathrm E$ ($\gamma \neq 0$). Therefore, $R^*(\mathrm E)=\lambda \mathrm E$ and 
$$
R(\mathrm{E})+R^*(\mathrm{E})=2\lambda \mathrm{E}+\theta e_{12},\ \theta\neq 0.
$$
Thus, $R$ doesn't satisfy \eqref{usE}. Using similar arguments, we get that operators from lines  2,3,4 do not satisfy \eqref{usE}.

Consider the line 7. In this case, $R^*(\mathrm{E})=\lambda \mathrm{E}+2\theta e_{21}$. Thus,
$$
R(\mathrm{E})+R^*(\mathrm{E})=2\lambda \mathrm{E}+2\theta e_{21}\neq \gamma \mathrm{E}.
$$

Similar arguments can be used to show that operators that are conjugate to operators from lines 8, 9,11,12 do not satisfy \eqref{usE}.

It is left to consider lines 5, 6, 10, 13. 

Suppose that $R$ is conjugate to an operator from  line 5. Then we have that $R^*(\mathrm{E})=\lambda \mathrm{E}+\alpha_1 h$. Thus, $R(\mathrm{E})+R^*(\mathrm{E})=\gamma \mathrm{E}$ if and only if $\alpha_1=-\theta$.

In 6, 10 and 13 it is easy to see that $R^*(\mathrm{E})=\lambda \mathrm{E}$. Thus, in this case, $R$ satisfies \eqref{usE} without a restriction. The proposition is proved.

From Proposition 3, it follows that it is enough to consider the following operators:

\begin{gather*}
 1.\  R(\mathrm{E})=\lambda \mathrm{E}+\theta h,\ R(h)=-\theta\mathrm{E}+\alpha_2h,\ R(e_{12})=-e_{12},\ R(e_{21})=0,\ \lambda,\alpha_2\in\mathbb C \ \theta\neq 0;\\
2.\  R(\mathrm{E})=\lambda \mathrm{E},\ R(h)=0,\ R(e_{12})=-e_{12}+th;\ R(e_{21})=0,\ t\in\{0,1\};\\
3.\  R(\mathrm{E})=\lambda \mathrm{E},\ R(h)=th,\ R(e_{21})=0,\ R(e_{12})=-e_{12},\ t\in \mathbb C,\ t\neq 0;\\
 4.\ R(\mathrm{E})=\lambda\mathrm{E}, R(h)=th,\ R(e_{12})=te_{12},\ R(e_{21})=te_{21},\ t\in\{0,-1\}.
\end{gather*}

We will consider operators 1-4 one by one.

{\bf Proposition 4.} Let $R$ be a Rota-Baxter operator of weight 1 defined as:
\begin{gather*}
   R(\mathrm{E})=\lambda \mathrm{E}+\theta h,\ R(h)=-\theta\mathrm{E}+\alpha_2h,\ R(e_{12})=-e_{12},\ R(e_{21})=0,\ \lambda,\alpha_2\in\mathbb C \ \theta\neq 0;
\end{gather*}
Then $R$ satisfy \eqref{rb3} if and only if $\alpha_2=-\frac{1}{2}$. In this case, for every $a\in sl_2(\mathbb C)$ we have $R(a)+R^*(a)+a=0$. Therefore, if $\alpha_2=-\frac{1}{2}$, then $R$ also satisfies \eqref{rb1}.

{\bf Proof.} 
Direct computations show that $R^*(\mathrm{E})=\lambda \mathrm{E}-\theta h$, $R(h)=\theta \mathrm{E}+\alpha_2 h$, $R^*(x)=0$ and $R^*(y)=-y$. 

Suppose that $R$ satisfies \eqref{rb3}. We have
$$
R([h,x])+R^*([h,x])=2R(x)+2R^*(x)=-2x.
$$
On the other hand, 
$$
[R(h),x]+[R^*(h),x]=2\alpha_2 x+2\alpha_2 x=4\alpha_2 x.
$$
 Therefore, $\alpha_2=-\frac{1}{2}$.
 
 Conversely, let $\alpha_2=-\frac{1}{2}$. It is easy to see that in this case, for all $a\in sl_2(\mathbb C)$ we have that
 $$
 R(a)+R^*(a)+a=0.
 $$
This means that equations \eqref{rb3} and \eqref{rb1} are true for $a,b\in sl_2(\mathbb C)$. Since $R(\mathrm{E})+R^*(\mathrm{E})\in Z(gl_2(\mathbb C))$, it is follows that equations \eqref{rb3} and \eqref{rb1} are true for all $a,b\in gl_2(\mathbb C)$. The proposition is proved.

{\bf Proposition 5.} Let $R$ be a Rota-Baxter operator of weight 1 defined as:
\begin{gather*}
 R(\mathrm{E})=\lambda \mathrm{E},\ R(h)=0,\ R(e_{12})=-e_{12}+th;\ R(e_{21})=0,\ t\in\{0,1\};
\end{gather*}
Then $R$ does not satisfy \eqref{rb3} for any  $t$ and $\lambda$.

{\bf Proof.}
For $R$ we have:
$$
R^*(\mathrm{E})=\lambda \mathrm{E},\ R^*(h)=te_{21},\ R^*(e_{12})=0,\ R^*(e_{21})=-e_{21},\ t\in\{0,1\}.
$$
Then,
$$
[R([h,e_{12}]+R^*([h,e_{12}])=-2e_{12}. 
$$
On the other hand,
$$
[R(h),e_{12}]+[R^*(h),e_{12}]=-th.
$$
Thus, $R$ does not satisfy \eqref{rb3}.

{\bf Proposition 6.} Let $R$ be a Rota-Baxter operator of weight 1 defined as:
\begin{gather*}
 R(\mathrm{E})=\lambda \mathrm{E},\ R(h)=th,\ R(e_{21})=0,\ R(e_{12})=-e_{12},\ t\in \mathbb C,\ t\neq 0;
\end{gather*}
Then $R$ satisfy \eqref{rb3} if and only if $t=-\frac{1}{2}$. In this case, for every $a\in sl_2(\mathbb C)$ we have $R(a)+R^*(a)+a=0$. Therefore, if $t=-\frac{1}{2}$, then $R$ also satisfies \eqref{rb1}.

{\bf Proof.} The proof is similar to the proof of Proposition 4. 

{\bf Proposition 7.} Let $R$ be a Rota-Baxter operator of weight 1 defined as:
\begin{gather*}
 \ R(\mathrm{E})=\lambda\mathrm{E}, R(h)=th,\ R(e_{12})=te_{12},\ R(e_{21})=te_{21},\ t\in\{0,-1\}.
\end{gather*}
Then for any $t\in \{0,-1\}$, $R$ satisfies \eqref{rb3}. Moreover, $R$ satisfies \eqref{rb1} if and only if $t=0$.

{\bf Proof.}
The first statement is obvious since the restriction of $R$ on $sl_2(\mathbb C)$ is equal to $t\cdot \mathrm{id}$.

If $t=0$, then $R(sl_2(\mathbb C))=0$. Since $[gl_2(\mathbb C),gl_2(\mathbb C)]=sl_2(\mathbb C)$, $R$ satisfies \eqref{rb1}.

If $t=-1$, then direct computations show that
$$
R(R([h,e_{12}]+R^*([h,e_{12}])+[h,e_{12}])=-2R(e_{12})\neq 0.
$$
Thus, if $t=-1$, then $R$ does not satisfy \eqref{rb1}.

Now we are ready to prove the main result of the paper:

{\bf Theorem 5.} Let $r=\sum a_i\otimes b_i\in gl_2(\mathbb C)\otimes gl_2(\mathbb C)$ be a solution of CYBE such that $r+\tau(r)\neq 0$ is $gl_2(\mathbb C)$-invariant. Then, up to the action of $Aut(gl_2(\mathbb C))$ and multiplication by a nonzero scalar, $r$ is equal to the one of the following:
\begin{gather}
  \label{r1} 1.\ r=\mathrm{E}\otimes(\lambda\mathrm{E}+\theta h)- h\otimes (\theta \mathrm{E}+\frac{1}{4} h)-e_{21}\otimes e_{12},\ \lambda\in \{0,1\},\ \theta\in \mathbb C, \theta\neq 0;\\
   2.\ r= \lambda \mathrm{E}\otimes \mathrm{E}-\frac{1}{2} h\otimes h-e_{21}\otimes e_{12},\ \lambda\in \{0,1\}; \\
  3.\ r= \lambda \mathrm{E}\otimes \mathrm{E},\ \lambda\in \{0,1\}. 
\end{gather}
    {\bf Proof.} Let $r\in gl_2(\mathbb C)\otimes gl_2(\mathbb C)$ be a solution of CYBE such that $r+\tau(r)$ is $gl_2(\mathbb C)$-invariant. By Theorem 2, the corresponding map $R$ defined as \eqref{rb2} is a Rota-Baxter operator of a nonzero weight. Up to a multiplication by a nonzero scalar we may assume that the weight of $R$ is equal to 1. Thus, $R$ satisfies \eqref{rb3} and \eqref{rb1}. From Propositions 3-6 it is follows that up to a conjugation with automorphisms from $Aut(M_2(\mathbb C))$, $R$ is one of the following:
    \begin{gather*}
 1.\  R(\mathrm{E})=\lambda \mathrm{E}+\theta h,\ R(h)=-\theta\mathrm{E}-\frac{1}{2}h,\ R(e_{12})=-e_{12},\ R(e_{21})=0,\ \lambda,\theta\in\mathbb C \ \theta\neq 0;\\
2.\ R(\mathrm{E})=\lambda \mathrm{E},\ R(h)=-\frac{1}{2}h,\ R(e_{21})=0,\ R(e_{12})=-e_{12},\ \lambda\in \mathbb C;\\
3. \ R(\mathrm{E})=\lambda\mathrm{E}, R(h)=R(e_{12})=R(e_{21})=0\ \lambda\in\mathbb C.
\end{gather*}

The second and the third maps do not depend on $\theta$, that means that they lie in different orbits modulo $Aut(M_2(\mathbb C))$. A simple check shows that maps of type 1  with different scalars $\theta_1$ and $\theta_2$ are not conjugate by means of $Aut(M_2(\mathbb C))$.

Therefore, up to the action of $Aut(M_2(\mathbb C))$ and multiplication by a nonzero scalar, $r$ is one of the following:
   \begin{gather*}
 1.\  \frac{1}{2}\mathrm{E}\otimes (\lambda \mathrm{E}+\theta h)+\frac{1}{2} h\otimes (-\theta\mathrm{E}-\frac{1}{2}h)-e_{21}\otimes e_{12},\ \lambda,\theta\in\mathbb C \ \theta\neq 0;\\
2.\ \lambda\mathrm{E}\otimes\mathrm{E}-\frac{1}{4}h\otimes h- e_{21}\otimes e_{12},\ \lambda\in \mathbb C;\\
3. \ \lambda \mathrm{E}\otimes\mathrm{E}, \lambda\in\mathbb C.
\end{gather*}

By Proposition 2, it is left to consider the action of an automorphism $\psi_{\beta}$ defined in \eqref{aut} for $\beta\in \mathbb C$, $\beta\neq 0$. Let $r$ be a solution of type 1 defined above. If $\lambda=0$, then we obtain the solution \eqref{r1} with $\lambda=0$. If $\lambda\neq 0$, then take $\beta=\frac{\sqrt{2}}{\sqrt{\lambda}}$
\begin{gather*}
    (\psi_{\beta}\otimes \psi_{\beta})r=\frac{1}{2}\beta \mathrm{E}\otimes (\lambda \beta \mathrm{E}+\theta h)+\frac{1}{2} h\otimes (-\theta\beta\mathrm{E}-\frac{1}{2}h)-e_{21}\otimes e_{12}\\
    =\mathrm{E}\otimes (\frac{\lambda \beta^2}{2} \mathrm{E}+\frac{\theta\beta}{2} h)- h\otimes (\frac{\theta\beta}{2}\mathrm{E}+\frac{1}{4}h)-e_{21}\otimes e_{12}\\
    =\mathrm{E}\otimes(\mathrm{E}+\theta_1 h)- h\otimes (\theta_1 \mathrm{E}+\frac{1}{4} h)-e_{21}\otimes e_{12}.
\end{gather*}
where $\theta_1=\frac{\theta\beta}{2}$. Thus, we obtain the solution \eqref{r1} for $\lambda=1$.

The rest of the cases can be considered similarly. The theorem is proved.

{\bf Corollary.} As a corollary of Theorem 5 we obtain a well known description of solutions $r$ of CYBE on $sl_2(\mathbb C)$ such that $r+\tau(r)$ is $sl_2(\mathbb C)$-invariant: up to a multiplication to a nonzero scalar and the action of $Aut(sl_2(\mathbb C))$, there is only one solution:
$$
r=\frac{1}{4}h\otimes h+e_{12}\otimes e_{21}.
$$

\end{document}